\documentclass[12pt]{article}
\usepackage{amssymb}
\usepackage{graphicx}
\usepackage{subcaption}
\usepackage{amsfonts}
\usepackage{latexsym}
\usepackage{amsthm}
\usepackage{amsmath}

\usepackage{amssymb, amscd}

\usepackage{url}
\usepackage[T1]{fontenc}
\usepackage{color}
\usepackage{tabu}

\usepackage{bbm}

\usepackage{amsmath,amssymb,amscd,}
\usepackage[all]{xy}
\usepackage{color}

\numberwithin{equation}{section}

\newtheorem{conj}{Conjecture}[section]
\newtheorem{thm}[conj]{Theorem}
\newtheorem{rem}[conj]{Remark}

\newtheorem{defin}[conj]{Definition}
\newtheorem{prop}[conj]{Proposition}
\newtheorem{cor}[conj]{Corollary}
\newtheorem{lema}[conj]{Lemma}

\begin{document}

 \title{Colored Tverberg problem,\\ extensions and new results}

\author{{Du\v{s}ko Joji\'{c}} \\
 {\small Faculty of Science}\\[-2mm]
 {\small University of Banja Luka}
\and Gaiane Panina\\
{\small St. Petersburg State University} \\[-2mm]
{\small St. Petersburg Department of}\\[-2mm] {\small Steklov Mathematical Institute}
\and Rade T. \v{Z}ivaljevi\'{c}\\
{\small Mathematical Institute}\\[-2mm] {\small SASA,
  Belgrade}\\[-2mm]}

\maketitle
\begin{abstract}\noindent
We prove  a ``multiple colored Tverberg theorem''  and  a ``balanced colored Tverberg theorem'', by applying different methods, tools and ideas.  The proof of the first theorem uses multiple chessboard complexes  (as configuration spaces) and   Eilenberg-Krasnoselskii theory of degrees of equivariant maps for non-free actions. The  proof of the second result relies on high connectivity of the configuration space, established by   discrete Morse theory.
\end{abstract}

\section{Introduction}
\label{sec:intro}

Tverberg-Van Kampen-Flores type results have been for decades one of the
central research themes in topological combinatorics. The last decade has been
particularly fruitful, bringing the resolution (in the negative) of the general ``Topological Tverberg Problem'' \cite{MaWa14, F,  bfz2, MaWa15, MaWa16}, as summarized  by several review papers \cite{BBZ-50, BS-50, Sk-18, Z17}.

New positive results include the proof \cite[Theorem 1.2]{jvz2} of the ``Balanced Van Kampen-Flores theorem'' and the development of ``collectively unavoidable simplicial complexes'', leading to very general theorems of Van Kampen-Flores type \cite{jpz1}.

\medskip
Somewhat surprisingly   the {\em colored Tverberg problem}, which in the past also occupied one of the central places  \cite{M, Z17}, doesn't seem to have been directly affected by these developments.

\medskip
For example the {\em Topological type A colored Tverberg theorem} (Theorem 2.2 in \cite{BMZ})
is known to hold if $r$ is a prime number (see Section \ref{sec:intro-AB} for the definition of the parameter $r$) and at present we don't know  what happens in other cases.
For comparison, the topological Tverberg theorem is known to be true if $r=p^\nu$ is a prime power, however we know today that this condition is essential (not an artefact of the topological methods used in the proof).

\medskip
Moreover, while there has been a notable, more recent  activity \cite{jvz2, jpz1} in the area of ``monochromatic'' Tverberg-Van Kampen-Flores theorems, the corresponding  new positive colored-type results seem to be virtually non-existent after the {\em arXiv} announcement of  \cite{BMZ}(more than ten years ago), where the type A colored Tverberg theorem was  established.

\bigskip

In this paper we open two new directions of studying colored Tverberg problem by proving a ``multiple Colored Tverberg theorem'' (Theorem \ref{thm:seven-points}) and ``balanced Colored Tverberg theorem'' (Theorem \ref{ThmColBalanced}).

\medskip
The ``multiple Colored Tverberg theorem'' has evolved from the simplified  proof of the type A colored Tverberg theorem given in \cite{vz11} and relies on the Eilenberg-Krasnoselskii theory of degrees of equivariant maps for non-free action.

\medskip
The ``balanced colored Tverberg theorem'' is a relative of the type B colored Tverberg theorem \cite{ZV92, VZ94} and the ``balanced Van Kampen-Flores theorem'' \cite{jvz2}.  It uses the methods of Forman's discrete Morse theory, and builds on the methods and ideas developed in our earlier papers \cite{jnpz, jpvz, jpz1}.

\bigskip\noindent
\textbf{Acknowledgements.} Proposition \ref{T:connect} is supported by the Russian Science Foundation under grant 16-11-10039. R. \v Zivaljevi\' c was supported by
the Ministry of Education, Science and Technological Development of Serbia (via grants to Mathematical Institute SASA).

\subsection{Brief overview of Tverberg's theorem and its relatives}
\label{sec:intro-AB}

``Tverberg type theorems'' is a common name for a growing family of
theorems, conjectures and problems about special partitions (patterns) of finite sets of points (point clouds) in the affine euclidean space  $\mathbb{R}^{d}$.

\medskip
The original Tverberg's theorem
 claims that every set $S\subset \mathbb{R}^d$
with $(r-1)(d+1)+1$ elements can be partitioned $S =
S_1\cup\ldots\cup S_r$ into $r$ nonempty, pairwise disjoint
subsets $S_1, \ldots, S_r$ such that the corresponding convex
hulls have a nonempty intersection:
\begin{equation}
\bigcap_{i=1}^r {\rm conv}(S_i)\neq\emptyset  \, .
\end{equation}

Among the numerous relatives (refinements, predecessors, generalizations) of Tverberg's theorem are the classical  Radon's lemma, Topological Tverberg theorem,  Colored Tverberg theorems, as well as the classical Van Kampen-Flores theorem  and its generalizations.

\medskip The original Tverberg's theorem  can be reformulated as the statement that for each linear (affine) map $f : \Delta^N \stackrel{a}{\longrightarrow}
\mathbb{R}^d$ from a $N$-dimensional simplex ($N =(r-1)(d+1)$) there exist $r$ nonempty disjoint
faces $\Delta_1,\ldots, \Delta_r$ of the simplex $\Delta^N$ such
that $f(\Delta_1)\cap\ldots\cap f(\Delta_r)\neq\emptyset$. This
form of Tverberg's result can be abbreviated as follows,
\begin{equation}\label{eqn:summary-Tverberg}
(\Delta^{(r-1)(d+1)} \stackrel{a}{{\longrightarrow}} \mathbb{R}^d)
\Rightarrow (r- {\rm intersection} )
\end{equation}
where the input on the left is (as indicated) an affine map and as the output on the right we  obtain the existence of a  ``Tverberg $r$-intersection'', that is a collection of $r$ disjoint faces of the simplex $\Delta^{(r-1)(d+1)}$ which overlap in the image.

\medskip
The reformulation (\ref{eqn:summary-Tverberg}) is very useful since it places the theorem in a broader context and motivates (potential) extensions and generalizations of Tverberg's (affine) theorem.

\medskip
For example the affine input map can be replaced by an arbitrary continuous map$f : \Delta^N\rightarrow \mathbb{R}^d$.
The corresponding more general statement (known as the Topological Tverberg theorem)
\begin{equation}\label{eqn:topological-Tverberg}
(\Delta^{(r-1)(d+1)} \longrightarrow \mathbb{R}^d)
\Rightarrow (r- {\rm intersection} )
\end{equation}
is also true  provided $r = p^\nu$ is a prime power.

\medskip
Another possibility is to prescribe in advance which pairwise disjoint faces $\Delta_1,\dots, \Delta_r$ of $\Delta^N$ are acceptable, preferred or ``admissible'', say by demanding that $\Delta_i\in K$ for a chosen simplicial subcomplex $K\subseteq \Delta^N$. The input map is now a continuous (affine, simplicial) map $f : K \rightarrow \mathbb{R}^d$ and the conclusion is the same as in the Topological Tverberg theorem.

\medskip
The following four statements are illustrative for results of
`colored Tverberg type' (see \cite{Z17} for more detailed presentation).

\begin{equation}\label{eqn:(1)}
(K_{3,3} \longrightarrow {\mathbb R}^2) \Rightarrow (2- {\rm
intersection})
\end{equation}
%\vspace{-7mm}
\begin{equation}\label{eqn:(2)}
(K_{3,3,3} \stackrel{a}{{\longrightarrow}} {\mathbb R}^2)
\Rightarrow (3- {\rm intersection})
\end{equation}
\begin{equation}\label{eqn:(3)}
(K_{5,5,5} \longrightarrow {\mathbb R}^3) \Rightarrow (3- {\rm
intersection})
\end{equation}
\begin{equation}\label{eqn:(4)}
(K_{4,4,4,4} \longrightarrow {\mathbb R}^3) \Rightarrow (4- {\rm
intersection})
\end{equation}

\medskip\noindent The complex $K_{t_0,t_1,\ldots,
t_k}=[t_0]\ast[t_1]\ast\ldots\ast[t_k]$ is by definition the
complete multipartite simplicial complex obtained as a join of
$0$-dimensional complexes (finite sets). For example $K_{p,q}=[p]\ast[q]$ is the complete bipartite graph obtained by connecting each of $p$ `red vertices'
with each of $q$ `blue vertices'.

\medskip
More generally a coloring of vertices of a simplex by $k+1$ colors is a partition  $V = Vert(\Delta^N) =  C_0\uplus C_1\uplus  \dots \uplus  C_k$  into ``monochromatic'' subsets $C_i$. A subset $\Delta\subseteq V$  is called a multicolored set or a {\em rainbow simplex} if and only if $\vert \Delta\cap C_i\vert \leq 1$ for each $i=0,\dots, k$. If the cardinality of $C_i$ is $t_i$ we observe that $K_{t_0,t_1,\ldots, t_k}$ is precisely the subcomplex of all rainbow simplices in $\Delta^N$.

\medskip
We refer the reader to \cite{ZV92, VZ94, BMZ}, \cite{M, Ziv-96-98} and       \cite{BBZ-50, BS-50, Sk-18, Z17}
for more general statements,  proofs, history and applications of colored Tverberg theorems.

\medskip
Following the classification proposed in \cite{Z17}  we say that a colored Tverberg  theorem is of {\em type A} if $k\geq d$ (where $k+1$ is the number of colors and $d$ is the dimension of the ambient euclidean space). In the case of the opposite inequality $k< d$ we say that it is of {\em type B}. The main difference between these two types is that in the type B case the number $r$ of intersecting rainbow simplices must satisfy the inequality $r\leq d/(d-k)$, while in the type A case there are no {\em a priori}  constraints on these numbers.

\medskip
In agreement with this classification (\ref{eqn:(2)}) and  (\ref{eqn:(4)}) are classified as  {\em topological type A colored Tverberg theorems} while (\ref{eqn:(1)}) and (\ref{eqn:(3)}) are instances of {\em topological type B colored Tverberg theorem}.

\medskip
The following general results (Theorems \ref{A} and \ref{B}) are the main representatives of these two classes of colored Tverberg theorems. In particular (\ref{eqn:(1)}) (\ref{eqn:(3)}) and (\ref{eqn:(4)}) are their easy consequences.

\medskip\noindent
 {\bf Caveat:} Here and elsewhere in the paper we do not distinguish the $N$-dimensional (geometric) simplex  $\Delta^N$ from a (combinatorial)  simplex $\Delta_{[m]} = 2^{[m]}$ (abstract simplicial complex) spanned by $m$ vertices (if $m = N+1$). In agreement with this convention,  subsets $S\subset [m]$ are interpreted as simplices, faces of $\Delta_{[m]}$. For $S\subset [m]$  we have ${\rm dim}(S) = \vert S\vert -1$ where $\vert S\vert$ is the cardinality of $S$.

\begin{thm}{\rm (Type A) \cite{BMZ}} \label{A}
Let $r\geq 2$ be a prime, $d\geq 1$, and $N:=(r-1)(d+1)$. Let $\Delta^N$ be an $N$-dimensional simplex with a partition (coloring) of its vertex set into $d+2$ parts, $V  = [N+1] = C_0\uplus \dots \uplus C_d\uplus C_{d+1}$, with $\vert C_i\vert = r-1$ for $i\leq d$ and $\vert C_{d+1}\vert =1$. Then for any continuous map $f : \Delta^N \rightarrow \mathbb{R}^d$, there are $r$ disjoint rainbow simplices  $\Delta_1,\dots, \Delta_r$ of $\Delta^N$ satisfying
\[
f(\Delta_1) \cap \dots \cap f(\Delta_r) \neq\emptyset \, .
\]
\end{thm}

\begin{thm}{\rm (Type B) \cite{VZ94, Ziv-96-98}}\label{B}
Assume that $r=p^\nu$ is a prime power,  $d\geq 1$, and $k<  d$.  Let $[m] = C_0\uplus \dots \uplus C_k$ be a coloring partition of vertices of the simplex $\Delta_{[m]}$  where $m = (2r-1)(k+1)$ such that $r\leq d/(d-k)$ and $\vert C_i\vert  = 2r -1$ for each $i$.
Then for any continuous map $f : \Delta_{[m]} \rightarrow \mathbb{R}^d$, there are $r$ disjoint rainbow simplices  $\Delta_1,\dots, \Delta_r$ of $\Delta_{[m]}$ satisfying
\[
f(\Delta_1) \cap \dots \cap f(\Delta_r) \neq\emptyset \, .
\]
\end{thm}

\subsection{Multiple colored Tverberg theorem}

The implication (\ref{eqn:(2)}) (Section~\ref{sec:intro-AB}) is an instance of a
result of B\'{a}r\'{a}ny and Larman \cite{BL92}. It says that each
collection of nine points in the plane, evenly colored by three
colors, can be partitioned into three  `rainbow triangles' which have a common point.

\medskip
At present it is not known if the following non-linear (topological) version of (\ref{eqn:(2)}) is true or not:

\begin{equation}\label{eqn:(2)-bis}
(K_{3,3,3} \longrightarrow {\mathbb R}^2)
\Rightarrow (3- {\rm intersection}) \, .
\end{equation}
In other words we don't know if for each continuous map $f : K_{3,3,3}\rightarrow \mathbb{R}^2$ there exist three pairwise vertex disjoint simplices $\Delta_1, \Delta_2, \Delta_3$ in $K_{3,3,3}$ such that $f(\Delta_1)\cap  f(\Delta_2)  \cap  f(\Delta_3)  \neq\emptyset$.

The implication (\ref{eqn:(2)-bis}) clearly follows from the following stronger statement:
\begin{equation}\label{eqn:(2)-bis-bis}
(K_{3,3,3,1} \longrightarrow {\mathbb R}^2)
\Rightarrow (4- {\rm intersection}) \, .
\end{equation}
However the implication (\ref{eqn:(2)-bis-bis}) is also not known to hold in full generality (and we strongly suspect that it is not the case).

\medskip
The following ``multiple Colored Tverberg theorem'' claims that the implication (\ref{eqn:(2)-bis-bis}) is true for continuous maps  $f : K_{3,3,3,1} \rightarrow {\mathbb R}^2$ which satisfy an additional (3-to-2) constraint. (The reader may find it instructive to  read first its affine version, Corollary~\ref{cor:seven-points}.)

\begin{thm}\label{thm:seven-points}
Suppose that $K = K_{3,3,3,1}\cong [3]\ast [3]\ast [3] \ast [1]$ is a $3$-dimensional simplicial complex on a ten-element vertex set, divided into four color classes.
  Assume that $f$ is a (3-to-2) map, meaning that $f = \widehat{f} \circ \alpha$ for some map $\widehat{f} : K_{2,2,2,1} \longrightarrow \mathbb{R}^2$ where
  \[
              \alpha : K_{3,3,3,1} \longrightarrow K_{2,2,2,1}
  \]
   is the simplicial map arising from a choice of epimorphisms $[3] \rightarrow [2]$. Then there exist four pairwise vertex disjoint simplices ($4$ vertex-disjoint rainbow simplices) $\Delta_1,\Delta_2,\Delta_3,\Delta_4$ in $K$ such that
\begin{equation}\label{eqn:4-intersection}
          f(\Delta_1)\cap f(\Delta_2)\cap f(\Delta_3)\cap f(\Delta_4)\neq\emptyset \, .
\end{equation}
\end{thm}

\begin{cor}\label{cor:seven-points}
Suppose that $X$ is a collection of $7$ points in the plane $\mathbb{R}^2$. Moreover, assume that these points are colored by $4$ colors, meaning that there is a partition
$X = A \sqcup B \sqcup C \sqcup D$  into monochromatic sets, where $ A =\{a_1, a_2\}, B = \{ b_1, b_2\}, C = \{ c_1, c_2\}$ are $2$-element sets and $D =\{d\}$ is a singleton.

Then there exist four rainbow sets $\Delta_i \subset A\cup B\cup C \cup D$ $(i=1,2,3,4)$ such that

\begin{enumerate}
\item[{\rm (1)}]
    \[
        Conv(\Delta_1) \cap  Conv(\Delta_2) \cap  Conv(\Delta_3) \cap Conv(\Delta_4) \neq\emptyset
    \]
\item[{\rm (2)}] each $a_1, b_1, c_1, d$ appears as a vertex in exactly one of the sets $\Delta_i$ and each $a_2, b_2, c_2$ appears as a vertex in exactly two of the sets $\Delta_i$.
\end{enumerate}
\end{cor}

The following corollary says that the implication (\ref{eqn:(2)-bis}) is true for a special class of non-linear maps.
\begin{cor}\label{cor:3-to-2}
Assume that $f : K_{3,3,3} \longrightarrow  \mathbb{R}^2$ is a continuous map which admits a factorization
\begin{equation}\label{eqn:factorization-1}
K_{3,3,3} \stackrel{\alpha}{\longrightarrow} K_{2,2,2} \stackrel{\widehat{f}}{\longrightarrow}  \mathbb{R}^2
\end{equation}
for some $\widehat{f}$, where $\alpha$  is a \mbox{ {\rm (3-to-2)}} map. Then there exist three disjoint triangles $\Delta_1, \Delta_2, \Delta_3$ in $K_{3,3,3}$ such that
\[
f(\Delta_1)\cap f(\Delta_2) \cap f(\Delta_3) \neq\emptyset \, .
\]
\end{cor}

 The proof of both Theorem~\ref{thm:seven-points} and its corollaries is postponed for Section~\ref{sec:rade}. The proofs rely on Eilenberg-Krasnoselskii theory of degrees of equivariant maps for non-free actions, see the monograph   \cite{KB} for a detailed presentation of the theory.

\subsection{Balanced colored Tverberg-type theorem}

Our ``balanced colored Tverberg theorem'' (Theorem \ref{ThmColBalanced})  was originally envisaged as an extension of the type B colored Tverberg theorem (Theorem \ref{B}) in the direction of the following theorem which is often referred to as the {\em balanced extension of the generalized Van Kampen-Flores theorem}.

\begin{thm}\label{thm:glavna-jvz} {\rm (\cite[Theorem 1.2]{jvz2})}
Let $r\geq 2$ be a prime power, $d \geq 1$, $N \geq (r - 1)(d +
2)$, and $rk+s \geq (r-1)d$ for integers $k \geq 0$ and $0 \leqslant s
< r$. Then for every continuous map $f : \Delta^N \rightarrow
\mathbb{R}^d$, there are $r$ pairwise disjoint faces
$\Delta_1,\ldots,\Delta_r$ of $\Delta^N$ such that
$f(\Delta_1)\cap \cdots\cap f(\Delta_r) \neq \emptyset$, with
$\textrm{dim }\Delta_i\leqslant k+1 $ for $1 \leqslant i \leqslant s$ and
$\textrm{dim }\Delta_i\leqslant k$ for $s < i \leqslant r$.
\end{thm}

If  one assumes that $(r-1)d$ is divisible by $r$, in which case $s=0$ and $\textrm{dim }\Delta_i\leqslant k$ for each $i$, then  Theorem \ref{thm:glavna-jvz} reduces to the `equicardinal'  generalized Van Kampen-Flores theorem of Sarkaria \cite{sar}, Volovikov \cite{Vol96-2} and Blagojevi\' c, Frick and Ziegler \cite{bfz1}.

\medskip
The following ``balanced colored Tverberg theorem'' can be described as a relative of Theorem \ref{thm:glavna-jvz} and ``balanced'' extension of Theorem \ref{B}.

\begin{thm} \label{ThmColBalanced}
Assume that $r=p^\nu$ is a prime power and let  $d\geq 1$. Let integers $k \geq 0$ and $0< s \leqslant  r$ be such that    \vspace{-0.4cm}
\begin{align}\label{eqn:ceiling}
r(k-1)+s &= (r-1)d,  \hbox{  or more explicitly,} \\
  k  := \lceil {(r-1)d}/{r}\rceil\quad &\mbox{{and} }\quad  s:=  (r-1)d-r(k-1) \, .
\end{align}
 Let $[m] = C_1\uplus \dots \uplus C_{k+1}$ be a coloring partition of vertices of  $\Delta_{[m]}$,  where $m = (2r-1)(k+1)$ and $\vert C_i\vert  = 2r -1$ for each $i$.
Then for any continuous map $f : \Delta_{[m]} \rightarrow \mathbb{R}^d$  there are $r$ disjoint rainbow simplices  $\Delta_1,\dots, \Delta_r$ of $\Delta_{[m]}$ satisfying $f(\Delta_1) \cap \dots \cap f(\Delta_r) \neq\emptyset$ such that
\begin{equation}\label{eqn:}
  \textrm{dim}(\Delta_i) = \vert \Delta_i\vert -1 \leqslant k  \mbox{ {\rm for} } 1 \leqslant i \leqslant s  \mbox{ \, {\rm and }\,  }  \textrm{dim}(\Delta_i)  \leqslant k - 1  \mbox{ {\rm for} } s < i \leqslant r \, .
\end{equation}
\end{thm}

\begin{rem}{\rm  Theorem \ref{B} is a special case of Theorem \ref{ThmColBalanced} for $s=r$. Indeed, the condition $r\leq d/(d-k)$  (in Theorem \ref{B}) is easily checked to be equivalent to the condition $rk \geq (r-1)d$.
  }
\end{rem}

 The proof of Theorem \ref{ThmColBalanced}  is based on high connectivity of the appropriate configuration space (Proposition \ref{T:connect}), which is proved by the methods of discrete Morse theory.

For the reader's convenience we briefly outline basic facts and ideas from discrete Morse theory in  Appendix 1. For a more complete presentation the reader is referred to  \cite{Forman2}.

 \section{Proof of the multiple Colored Tverberg theorem}
 \label{sec:rade}

 The first step in the proof of the multiple Colored Tverberg theorem (Theorem \ref{thm:seven-points}) is a standard reduction, via the {\em Configuration Space/Test Map scheme} \cite{Z17, M, Ziv-96-98}, to a problem of equivariant topology.

   \medskip Starting with a continuous map $f : K_{3,3,3,1} \rightarrow \mathbb{R}^2$ we build the associated configuration space as the deleted join
   $$
   (K_{3,3,3,1})^{\ast 4}_\Delta  =  ([3]\ast [3]\ast [3] \ast [1])^{\ast 4}_\Delta \cong
   (\Delta_{3,4})^{\ast 3}\ast [4]
   $$
 where $\Delta_{3,4}$ is the standard chessboard complex of all non-taking rook placements on a $(3\times 4)$-chessboard.

 \medskip
 The associated test map, designed to test if a simplex $\tau = (\Delta_1, \Delta_2, \Delta_3, \Delta_4)\in (K_{3,3,3,1})^{\ast 4}_\Delta$ satisfies the condition  (\ref{eqn:4-intersection}), is defined as a $\Sigma_4$-equivariant map
 \begin{equation}\label{eqn:test-map-Phi}
   \Phi :  (K_{3,3,3,1})^{\ast 4}_\Delta  \longrightarrow (\mathbb{R}^2)^{\ast 4}/D \hookrightarrow (W_4)^{\oplus 3}
 \end{equation}
where $D\subset (\mathbb{R}^2)^{\ast 4}$ is the diagonal ($2$-dimensional) subspace and $W_4$ is the standard, $3$-dimensional real permutation representation of $\Sigma_4$.

\medskip
Summarizing, in order to show that there exists a $4$-tuple $(\Delta_1, \Delta_2, \Delta_3, \Delta_4)$ satisfying (\ref{eqn:4-intersection}) it is sufficient to prove that the  $\Sigma_4$-equivariant map (\ref{eqn:test-map-Phi})  must have a zero.

\medskip
For the next step we need to use, side by side with  the standard chessboard complex $\Delta_{3,4}$ (on a $(3\times 4)$-chessboard),   a `multiple chessboard complex' $\Delta_{2,4}^{\mathbbm{1};\mathbb{L}})$, defined as the complex of all rook placements on a $(2\times 4)$-chessboard such that in the first column up to two rooks are permitted, while in all  rows and in the second column at most one rook is allowed.

 More general `multiple chessboard complexes' are studied in \cite{jvz}, and the notation follows this paper. In particular the vectors $\mathbbm{1} = (1,1,1,1)$ (respectively $\mathbb{L} = (1,2)$) describe the restriction on the rook placements in the rows (respectively the columns) of the $(2\times 4)$-chessboard.

\medskip
Now we use the fact that $f$  satisfies the (3-to-2) condition, which allows us to prove the  following lemma.

\begin{lema}\label{lema:factorization}
  Assume that $f$ is a (3-to-2) map, meaning that $f = \widehat{f} \circ \alpha$ for some map $\widehat{f} : K_{2,2,2,1} \longrightarrow \mathbb{R}^2$ where
  \[
              \alpha : K_{3,3,3,1} \longrightarrow K_{2,2,2,1}
  \]
   is the simplicial map arising from a choice of epimorphisms $[3] \rightarrow [2]$.        Under this condition the equivariant map (\ref{eqn:test-map-Phi}) admits a factorization     $\Phi  =  \widehat{\Phi} \circ \pi$  into $\Sigma_4$-equivariant maps,  as displayed in the following commutative diagram
  \begin{equation}\label{CD:3-to-2}
\begin{CD}
(\Delta_{2,4}^{\mathbbm{1};\mathbb{L}})^{\ast (3)}\ast [4] @>\widehat{\Phi}>> (W_4)^{\ast (3)}\\
@A\pi AA @A\cong AA\\
(\Delta_{3,4})^{\ast 3}\ast [4] @>\Phi>> (W_4)^{\ast (3)}
\end{CD}
\end{equation}
where $\Delta_{2,4}^{\mathbbm{1};\mathbb{L}}$ is the multiple chessboard complex defined above and $\pi$ is an epimorphism.
\end{lema}

\medskip\noindent
{\bf Proof:} The proof is by elementary inspection. Note that the map $\widehat{\pi} : \Delta_{3,4} \rightarrow \Delta_{2,4}^{\mathbbm{1};\mathbb{L}}$ , which induces the map $\pi$ in the diagram (\ref{CD:3-to-2}), is informally described as the map which unifies two columns of the $(3\times 4)$-chessboard into one column of the $(2\times 4)$-chessboard. \hfill $\square$

\medskip
Summarizing the first two steps we observe that the proof of Theorem \ref{thm:seven-points} will be complete if we show that the $\Sigma_4$-equivariant map $\widehat{\Phi}$ always has a zero. (Here we tacitly use the fact that $\pi$ is an epimorphism.)

\subsection{Equivariant maps from the manifold  $(\Delta_{2,4}^{\mathbbm{1};\mathbb{L}})^{\ast (3)}$}

The $\Sigma_4$-representation $W_4$ can be described as $\mathbb{R}^3$ with the action coming from the symmetries of regular tetrahedron $\Delta_{[4]}$, centered at the origin. If the map $\widehat{\Phi}$ has no zeros than there exists a $\Sigma_4$-equivariant map $$g : (\Delta_{2,4}^{\mathbbm{1};\mathbb{L}})^{\ast (3)}\ast [4]   \longrightarrow (\partial\Delta_{[4]})^{\ast (3)}$$
where $\partial\Delta_{[4]}$ is the boundary sphere of the simplex $\Delta_{[4]}$. This is ruled out by the following theorem.

\begin{thm}\label{thm:seven-BU}
 Let $G = (\mathbb{Z}_2)^2 = \{1, \alpha, \beta, \gamma\}$ be the {\em Klein four-group}. Let $\Delta_{2,4}^{\mathbbm{1};\mathbb{L}}$ be the multiple chessboard complex (based on a $2\times 4$ chessboard), where $ \mathbbm{1} = (1,1,1,1)$ and $\mathbb{L} = (2,1)$, and let $\partial\Delta_{[4]}\cong S^2$ be the boundary of a simplex spanned by vertices in $[4]$.
 Both $\Delta_{2,4}^{\mathbbm{1};\mathbb{L}}$ and $\partial\Delta_{[4]}\cong S^2$ are $G$-spaces, where the first action permutes the rows of the chessboard $[2]\times [4]$, while the second permutes the vertices of the $3$-simplex   $\Delta_{[4]}$. Under these conditions there does not exist a $G$-equivariant map
 \[
      f :   (\Delta_{2,4}^{\mathbbm{1};\mathbb{L}})^{\ast (3)}\ast [4] \longrightarrow (\partial\Delta_{[4]})^{\ast (3)}\cong (S^2)^{\ast (3)} \cong S^8
 \]
 where the joins have the diagonal $G$-action.
\end{thm}

Theorem \ref{thm:seven-BU} is proved by an argument involving the degree of equivariant maps which can be traced back to Eilenberg and Krasnoselskii, see \cite{KB} for a thorough  treatment and Appendix 1 for the statement of one of the main theorems.

\medskip
Before we commence the proof of  Theorem \ref{thm:seven-BU} let us describe a convenient geometric model for the complex $\Delta_{2,4}^{\mathbbm{1};\mathbb{L}}$. Recall that the Bier sphere $Bier(K)$ of a simplicial complex $K\subset 2^{[m]}$ is the deleted join $K\ast_\Delta K^\circ$ of $K$ and its Alexander dual $K^\circ$, see \cite{M} for more details on this subject.

\begin{lema}\label{lema:small-Bier}
The multiple chessboard complex $\Delta_{2,4}^{\mathbbm{1};\mathbb{L}}$ is a triangulation of a $2$-sphere. More explicitly, there is an isomorphism
$\Delta_{2,4}^{\mathbbm{1};\mathbb{L}} \cong Bier(\Delta_{[4]}^{(1)})$, where $\Delta_{[4]}^{(1)}$  is the $1$-skeleton of the tetrahedron $\Delta_{[4]}$ and
$Bier(K) = K\ast_\Delta K^\circ$ is the Bier sphere associated to a simplicial complex $K$ (and its Alexander dual $K^\circ$).
\end{lema}

\medskip The proof of Lemma~\ref{lema:small-Bier} is straightforward and relies on the observation that the subcomplexes of  $\Delta_{2,4}^{\mathbbm{1};\mathbb{L}}$,
generated by the vertices in the first (second) column of the chessboard $[2]\times [4]$, are respectively $K = \Delta_{[4]}^{(1)}$  and $K^\circ = (\Delta_{[4]}^{(1)})^\circ  = \Delta_{[4]}^{(0)}$.
The following lemma clarifies the structure of the sphere $\Delta_{2,4}^{\mathbbm{1};\mathbb{L}}$ as a $G$-space where $G = (\mathbb{Z}_2)^2 = \{1, \alpha, \beta, \gamma\}$  is the Klein four-group.

\begin{lema}\label{lema:G-Bier}
As a $G$-space the sphere $\Delta_{2,4}^{\mathbbm{1};\mathbb{L}}$ is homeomorphic to the regular octahedral sphere (positioned at the origin), where the generators $\alpha, \beta, \gamma$ are interpreted
as the $180^\circ$-rotations around the axes connecting pairs of opposite vertices of the octahedron.

More explicitly, let $\mathbb{R}^1_\alpha$ be the $1$-dimensional $G$-representation characterized by the conditions $\alpha x =  x, \beta x = \gamma x = -x $ ($\mathbb{R}^1_\beta$ and $\mathbb{R}^1_\gamma$
are defined similarly)   and let $S^0_\alpha, S^0_\beta, S^0_\gamma$ be the corresponding $0$-dimensional $G$-spheres. Then $\Delta_{2,4}^{\mathbbm{1};\mathbb{L}}$ is $G$-isomorphic to the $2$-sphere
$S(\mathbb{R}^1_\alpha \oplus \mathbb{R}^1_\beta \oplus \mathbb{R}^1_\gamma) \cong S^0_\alpha \ast S^0_\beta \ast S^0_\gamma$ with the induced $G$-action.
\end{lema}

\begin{rem}\label{rem:G-iso}
{\rm   Here is a geometric interpretation (visualization) of the  $G$-isomorphism $\Delta_{2,4}^{\mathbbm{1};\mathbb{L}} \cong Bier(\Delta_{[4]}^{(1)})$.  The geometric realizations of $K = \Delta_{[4]}^{(1)}$ and its Alexander dual $K^\circ = \Delta_{[4]}^{(0)}$ are respectively
constructed in the tetrahedron $\Delta_{[4]}$ and its polar body  $\Delta_{[4]}^\circ$.
If both tetrahedra are inscribed in the cube $I^3$, the geometric realization of $Bier(K)$ is naturally interpreted as a triangulation of
the boundary $\partial (I^3) $ of the cube $I^3$.
}
\end{rem}

\begin{lema}\label{lema:G-tetra}
As a $G$-space the boundary sphere of the tetrahedron $\partial\Delta_{[4]}$ is also homeomorphic to the octahedral sphere described in Lemma \ref{lema:G-Bier}. Moreover, there is a radial $G$-isomorphism $\rho :  \partial (I^3)\rightarrow  \partial\Delta_{[4]}$.
\end{lema}

Summarizing we conclude that the $G$-sphere we are studying in this section has two combinatorial   $\Delta_{2,4}^{\mathbbm{1};\mathbb{L}}, \partial\Delta_{[4]} = 2^{[4]}\setminus \{[4]\}$
and three equivalent geometric incarnations, the boundary of the cube $\partial  (I^3)$, the boundary of the tetrahedron  $\partial\Delta_{[4]}$, and the boundary of the octahedron $S^0_\alpha \ast S^0_\beta \ast S^0_\gamma$.

\subsection{Completion of the proof of Theorem \ref{thm:seven-BU}}

\begin{prop}\label{prop:deg-odd}
Let $\phi :  (\Delta_{2,4}^{\mathbbm{1};\mathbb{L}})^{\ast (3)} \rightarrow (\partial\Delta_{[4]})^{\ast (3)}$  be an arbitrary $G$-equivariant map. Then $${\rm deg}(\phi) \equiv 1 \, (\mbox{\rm mod}\, 2) \, .$$
\end{prop}

\medskip\noindent
{\bf Proof:} It follows from Theorem~\ref{thm:K-B-2.1} that ${\rm deg}(\phi) \equiv {\rm deg}(\psi) \, (\mbox{\rm mod}\, 2)$  for each two equivariant maps of the indicated spaces.  Here we use that fact that  $(\Delta_{2,4}^{\mathbbm{1};\mathbb{L}})^{\ast (3)} \cong  (S^2)^{\ast (3)} \cong S^8$ is a topological manifold.

Hence, it is enough to exhibit a single map $\psi$ with an odd degree. In light of the results from the previous section $(\Delta_{2,4}^{\mathbbm{1};\mathbb{L}})^{\ast (3)}$ and  $(\partial\Delta_{[4]})^{\ast (3)}$ are $G$-homeomorphic, $8$-dimensional spheres. If we choose the $G$-isomorphism $\psi : (\Delta_{2,4}^{\mathbbm{1};\mathbb{L}})^{\ast (3)} \rightarrow (\partial\Delta_{[4]})^{\ast (3)}$ then ${\rm deg}(\psi) = \pm 1$.
\hfill $\square$

\bigskip\noindent
{\bf Proof of Theorem \ref{thm:seven-BU}:}

\begin{equation}
\begin{CD}
(\Delta_{2,4}^{\mathbbm{1};\mathbb{L}})^{\ast (3)}\ast [4] @>f>> (\partial\Delta_{[4]})^{\ast (3)}\\
@AeAA @A\cong AA\\
(\Delta_{2,4}^{\mathbbm{1};\mathbb{L}})^{\ast (3)} @>\phi>> (\partial\Delta_{[4]})^{\ast (3)}
\end{CD}
\end{equation}
Suppose that a $G$-equivariant map $f$ exists. Let $e$ be the inclusion map and let $\phi = f\circ e$ be the composition.

The map $e$ is homotopically trivial, since ${\rm Image}(e) \subset {\rm Cone}(v)$ for each $v\in [4]$. This is, however, a contradiction since in light of Proposition~\ref{prop:deg-odd} the map $\phi$ has an odd degree.
\hfill $\square$

 \section{Proof of the balanced Color Tverberg theorem}
 \label{sec:dusko}

\medskip
Following the ``join variant'' of the {\em configuration space/test map}-scheme \cite{M} \cite{Z17},  a {\em configuration space} $\mathfrak{C} \subseteq \Delta_{[m]}^{\ast (r)}$, appropriate for the proof of Theorem \ref{ThmColBalanced}, collects together all joins $A_1\ast\dots\ast A_r:= A_1\uplus\dots\uplus A_r$ of disjoint rainbow simplices $A_i\subset [m]$, satisfying (after a permutation of indices) the condition   (\ref{eqn:}) from Theorem \ref{ThmColBalanced}.
For the future reference we record a more detailed  description of this configuration space.

\begin{defin}\label{def:conf-space}
The configuration space $\mathfrak{C}$ of $r$-tuples of disjoint rainbow simplices satisfying the  restrictions listed in Theorem \ref{ThmColBalanced} is the simplicial complex whose simplices are labeled by $$(A_1,...,A_r;B)$$
where\begin{itemize}
       \item $[m]=A_1\sqcup...\sqcup A_r\sqcup B$ is a partition such that $B\neq[m]$.
       \item Each  $A_i$ is a rainbow set (simplex), in particular  $|A_i|\leq k+1$ for each $i\in [r]$.
       \item The number of simplices $A_i$ with $|A_i|= k+1$ does not exceed $s$.
     \end{itemize}
\end{defin}

Note that the dimension of a simplex $(A_1,...,A_r;B)$  is  $|A_1|+...+|A_r|-1$. Moreover, a facet of a simplex $(A_1,...,A_r;B)$ is formally obtained by moving an element of some $A_i$ to $B$.

\begin{prop}\label{T:connect}
The configuration space $\mathfrak{C}$ is
$(rk+s-2)$-connected.
\end{prop}

Let us briefly explain how Theorem \ref{ThmColBalanced} can be deduced from  Proposition \ref{T:connect}. This is  a standard argument used for example in the proof of topological Tverberg theorem, see \cite[Section 6]{M} or \cite{Z17}.

\medskip
Suppose Theorem \ref{ThmColBalanced} is not true, which means  that  $f(A_1) \cap \dots \cap f(A_r) = \emptyset$ for all collections  $A_1,\dots, A_r$  of $r$ disjoint rainbow simplices  satisfying  (\ref{eqn:}). From here we deduce that there exists a $(\mathbb{Z}/p)^\nu$- equivariant mapping $$\Psi_f :\mathfrak{C}\rightarrow \mathbb{R}^{(d+1)r}$$  missing the diagonal
$D=\{(y,y,...,y): y\in \mathbb{R}^{d+1}\}$.

\medskip
However ${\rm Image}(\Psi_f)\subset \mathbb{R}^{(d+1)r}\setminus D$ contradicts Volovikov's theorem \cite{Vol96-2}, since \newline $\mathbb{R}^{(d+1)r}\setminus D$ is $(\mathbb{Z}/p)^\nu$-homotopy equivalent to a sphere of dimension $(r-1)(d+1)-1 = rk+s-2$ and the configuration space $\mathfrak{C}$ is by Proposition \ref{T:connect} $(rk+s-2)$-connected.

\bigskip

\medskip\noindent
{\bf Proof of Proposition \ref{T:connect}:}
Let us begin by introducing some useful abbreviations.

A set $A\subset [m]$  is called \textit{$C_i$-full} if it contains a vertex colored by $C_i$.

A simplex $(A_1,...,A_r;B)$  is called \textit{$C_i$-full} if each of the $A_i$ is full, or equivalently, if $|\bigcup_{i=1}^rA_i \bigcap C_i|=r$.

A simplex $(A_1,...,A_r;B)$  is $(k+1)$\textit{-full} if it contains (the maximal allowed  number) $s$  of $k+1$-sets among $A_i$.

A simplex $(A_1,...,A_r;B)$  is \textit{saturated} if it is $(k+1)$-full, and $|A_i|\geq k \ \ \forall i$.

 Saturated simplices are maximal faces of the configuration space  $\mathfrak{C}$.
Their dimension is $rk+s-1$.

\medskip

We now define a Morse matching for $\mathfrak{C}$. For a given simplex
$(A_1,...,A_r;B)$ we either describe a simplex  that is paired with it, or alternatively
recognize  $(A_1,...,A_r;B)$ as a critical simplex.

This is done stepwise. We shall have $r$ ``big'' steps, each of them further splitting  into consecutive $k+1$ small steps. Big steps treat  the sets $A_i$ one by one, and small steps treat colors one by one.

\medskip

\textbf{Step 1.}
\begin{description}
  \item[Step 1.1] Assume that the vertices of each color are enumerated by $\{1,2,...,2r-1\}$.
  Set $$a_1^1=min\Big[(A_1\cup B)\cap C_1\Big]$$
  and match $(A_1\cup a_1^1,A_2,...,A_r;B)$ with $(A_1,A_2,...,A_r;B\cup a_1^1)$  whenever both these simplices
  are elements of $\mathfrak{C}$.

  A simplex of type $(A_1\cup a_1^1,A_2,...,A_r;B)\in \mathfrak{C}$ is not matched iff it equals $$(\{a_1^1\},\emptyset,...,\emptyset;[m]\setminus  \{a_1^1\}).$$  It is  $0$-dimensional and it will stay unmatched until the end of the matching process.

   If a simplex of type $(A_1,A_2,...,A_r;B\cup a_1^1)$ is unmatched then
   either $A_1$ is $C_1$-full, or $|A_1|=k$, and $(A_1,A_2,...,A_r;B\cup a_1^1)$ is $(k+1)$-full.

  \item[Step 1.2]  Set $$a_1^2=min\Big[(A_1\cup B)\cap C_2\Big]$$
  and match $(A_1\cup a_1^2,A_2,...,A_r;B)$ with $(A_1,A_2,...,A_r;B\cup a_1^2)$  whenever both these are elements of $\mathfrak{C}$  that have not been matched on the Step 1.1.

  \begin{itemize}
    \item If a simplex of type $(A_1,A_2,...,A_r;B\cup a_1^2)$ is unmatched, then
   either $A_1$ is $C_2$-full, or $|A_1| = k$, and $(A_1,A_2,...,A_r;B\cup a_1^2)$ is $(k+1)$-full.

   Such simplices are called "\textit{Step 1.2-Type 1 unmatched} simplices".
    \item If a simplex of type $(A_1\cup a_1^2,A_2,...,A_r;B)$  is not matched, then $|A_1\cup a_1^2|=k+1$, and
  $(A_1\cup a_1^2,A_2,...,A_r;B)$ is $(k+1)$-full (these are  necessary but not sufficient conditions). The reason is that in this case  $(A_1,A_2,...,A_r;B\cup a_1^2)$  belongs to $\mathfrak{C}$  but might be matched on the Step 1.1.

Such simplices are called "\textit{Step 1.2-Type 2 unmatched} simplices."
  \end{itemize}

\medskip

In the sequel we use similar abbreviations. Step $i.j$ -- Type 1 means, that one cannot move an element colored by $j$ from $B$ to  $A_i$.
Step $i.j$ -- Type 2 means, that one cannot move an element colored by $j$ from
 $A_i$ to $B$.

  \item[Step 1.3] and subsequent steps (up to Step $1.k+1$)  go analogously.
\end{description}

Summarizing, we conclude:

  \begin{lema}\label{Lema Step1Unmatched}
     With the exception of the unique zero-dimensional unmatched simplex, if a simplex $(A_1,...,A_r;B)$ is unmatched after Step 1 then
       one of the following is valid:
     \begin{enumerate}
       \item either $|A_1|=k+1$, or
       \item $|A_1|=k$, and $(A_1,...,A_r;B)$ is $(k+1)$-full.
     \end{enumerate}
\end{lema}
Proof  follows directly from the  analysis of matching algorithm on small steps. \qed

\medskip

\medskip

\textbf{Step 2.}  Now we treat $A_2$  for the simplices that remained unmatched after Step 1.
\begin{description}
  \item[Step 2.1]
  Set $$a_2^1=min\Big[\Big((A_2\cup B)\setminus [1,a_1^1]\Big)\cap C_1\Big]$$
  and match $(A_1,A_2\cup a_2^1,...,A_r;B)$ with $(A_1,A_2,...,A_r;B\cup a_2^1)$  whenever both these are elements of $\mathfrak{C}$  that are not matched on Step 1.

\begin{itemize}
  \item If  a simplex $(A_1,A_2,...,A_r;B\cup a_2^1)$ is not matched now, then either $|A_2|=k$, and $(A_1,A_2,...,A_r;B\cup a_2^1)$ is $(k+1)$-full, or $A_2$ is $C_1$-full.

Such simplices will be called \textit{Step 2.1 --- Type 1 }simplices.
  \item
  If a simplex of type $(A_1,A_2\cup a_2^1,...,A_r;B)$ is not matched,  then it is $(k+1)$-full, and $|A_2|=k+1.$

Such simplices will be called \textit{Step 2.1 --- Type 2 }simplices.
\end{itemize}

  \item[Step 2.2] Set $$a_2^2=min\Big[\Big((A_2\cup B)\setminus [1,a_1^2]\Big)\cap C_2\Big]$$
  and match $(A_1,A_2\cup a_2^2,...,A_r;B)$ with $(A_1,A_2,...,A_r;B\cup a_2^2)$  whenever both these are elements of $\mathfrak{C}$  that are not matched before, that is, on Step 1, and on Step 2.1.

  \item[Step 2.3] and subsequent steps (up to Step $2.k+1$)  go analogously.
\end{description}

Summarizing, we conclude:

  \begin{lema}\label{LemStep2Unmatched}
     With the exception of the unique zero-dimensional unmatched simplex, if a simplex $(A_1,...,A_r;B)$ is unmatched after Step 2 then it is unmatched after Step 1 (for this we have Lemma \ref{Lema Step1Unmatched}), and also
       one of the following is valid:
     \begin{enumerate}
       \item either $|A_2|=k+1$, or
       \item $|A_2|=k$, and $(A_1,...,A_r;B)$ is $(k+1)$-full.\qed
     \end{enumerate}
\end{lema}

\medskip

 Steps 3,4,..., and  $r-1$ go analogously.

 \begin{lema}For all the steps $j=1,2,...,r-1$, the numbers $a_j^i$ are well-defined.
 \end{lema}
 Proof. Indeed, for $(A_1,...,A_r;B) \in \mathfrak{C}$, the set $B\cap C_i$ contains at least $r-1$ points.  (Here we use that $|C_i|=2r-1$ and $\vert A_j\cap C_i\vert \leq 1$ for each $j$.)
 The entries $a_1^i, a_2^i,...,a_{j-1}^i$  are either not in $B\cap C_i$, or (by construction) are the smallest consecutive entries of
 $B\cap C_i$. Altogether there are strictly less than $r-2$ of them.\qed

 \medskip

 A special attention should be paid to the last Step $r$.

 First, let us observe that (by construction) we already have:

 \begin{lema}
     With the exception of the unique zero-dimensional unmatched simplex, if a simplex $(A_1,...,A_r;B)$ is unmatched after Step $r-1$ then
       one of the following is valid:
     \begin{enumerate}
       \item  $|A_1|=|A_2|=...=|A_{r-1}|=k+1$, or
       \item  for some $i$, $|A_i|=k$, and $(A_1,...,A_r;B)$ is $(k+1)$-full.
     \end{enumerate}
\end{lema}

Proof: This follows from Lemma
\ref{LemStep2Unmatched} and its analogs for Steps $1,...,r-1$.\qed

\bigskip

\textbf{Step r.}  Now we turn our attention to  $A_r$.
\begin{description}
  \item[Step r.1]
  Set $$a_r^1=min\Big[\Big((A_r\cup B)\setminus [1,a_{r-1}^1]\Big)\cap C_1\Big].$$
  It might happen that the set $\Big[\Big((A_r\cup B)\setminus [1,a_{r-1}^1]\Big)\cap C_1\Big]$ is empty for $(A_1,...,A_r;B)$,
  so $a_r^1$  is ill-defined.

  This means that $(A_1,...,A_r;B)$ is $C_1$-full.
  Such simplex is left unmatched and called  \textit{Step $r.1$  --- Type 3} simplex.

If $a_r^1$ is well-defined, we proceed  in our standard way:
   we match $(A_1,A_2,...,A_r\cup a_r^1;B)$ with $(A_1,A_2,...,A_r;B\cup a_r^1)$  whenever both these are elements of $\mathfrak{C}$  that are not matched before.

  \item[Step r.2] Set $$a_r^2=min\Big[((A_r\cup B)\setminus [1,a_{r-1}^2])\cap C_2\Big].$$

  Again, if this number is ill-defined, this means that $(A_1,...,A_r;B)$ is $C_2$-full, and we leave the simplex Step $r.2$ --- Type 3 unmatched.

  Otherwise we proceed standardly.

  \item[Step r.3] and subsequent steps (up to Step r.k+1)  go analogously.
\end{description}

Summarizing, we conclude:

  \begin{lema}\label{LemStepUnmatched}
     With the exception of the unique zero-dimensional unmatched simplex, if a simplex $(A_1,...,A_r;B)$ is unmatched after Step $r$, then it is saturated.
\end{lema}
Proof.
Firstly, by Lemma \ref{LemStepUnmatched}, $|A_i|\geq k$ for $i=1,...,r-1$.

If a simplex $(A_1,...,A_r;B)$ has $|A_i|<k$ for some $i$, then it has some missing color. Let the smallest missing color be $j$. Then the simplex gets matched at Step $i.j$, since $a^i_j$ is well defined and can be added to $A_i$.

On each Step $i.j$, the simplex  $(A_1,...,A_r;B)$ is either Type 1, or Type 2, or (this might happen for Step $r.j$)
Type 3.  If it is at least once of Type 2 (does not matter on which step), then (by the same lemma) it is $(k+1)$-full, and therefore saturated.

If it is always of Type 1 on steps $1,...,r-1$ and not saturated, then $|A_i|=k+1$ for all $i=1,...,r-1$.

Since $s<r$, it is saturated.
\bigskip

It remains to prove the acyclicity of the matching.

Assume we have a gradient path

$$\alpha_0^p \nearrow \beta_0^{p+1} \searrow \alpha_1^p \nearrow \beta_1^{p+1} \searrow \alpha_2^p \nearrow \beta_2^{p+1} \searrow \,  \cdots \,   \searrow\alpha_m^p \nearrow \beta_m^{p+1} \searrow \alpha_{m+1}^p$$

For each of the simplices $\alpha$
consider the sequence of  numbers $$\Pi(\alpha):=(a_1^1,a_1^2,...,a_1^{k+1}, a_2^1, ...,a_2^{k+1},..., a_r^1, ...,a_r^{k+1})$$ These are all the numbers $a^i_j$ listed in the order that they appear in the matching algorithm;
 they  are well-defined including  the step where $\alpha$ gets matched.  If $a^i_r$ is ill-defined, set it  to be $\infty$.

\begin{lema} During a path, $\Pi(\alpha)$ strictly decreases w.r.t.\ lexicographic order. Therefore, the matching is acyclic.
\end{lema}
Firstly, it suffices to look at the two-step paths only:
$$\alpha_0^p \nearrow \beta_0^{p+1} \searrow \alpha_1^p \nearrow \beta_1^{p+1}$$
The proof is via an easy case analysis.  Here are two examples of how it goes:

(1) Assume  that $\alpha_0^p \nearrow \beta_0^{p+1}$ means adding color $i$ to $A_j$, and $\beta_0^{p+1} \searrow \alpha_1^p$ means removing color $i'>i$ from $A_j$.  Then
\begin{enumerate}
  \item either $\alpha_1^p$ is matched with some $p-1$-dimensional simplex
obtained by removing color $i$ from $A_j$,  and the path terminates right here,  or
  \item or $\alpha_1^p$ is matched before  Step $j.i$.
\end{enumerate}

(2) Assume  that $\alpha_0^p \nearrow \beta_0^{p+1}$ means adding color $i$ to $A_j$, and $\beta_0^{p+1} \searrow \alpha_1^p$ means removing color $i'$ from $A_{j'}$  with $j'<j$.  Then
\begin{enumerate}
  \item either $\alpha_1^p$ is matched   by adding  color $i'$ to $A_{j'}$,    or
  \item or $\alpha_1^p$ is matched before  Step $j'.i'$.
\end{enumerate}

\bigskip

This completes the proof of Theorem \ref{ThmColBalanced}.

\section*{Appendix 1. Discrete Morse theory}
\label{sec:DMT}

A discrete Morse function (or a discrete vector field) on a simplicial complex $K\subseteq 2^V$ is, by definition, an acyclic
matching on the Hasse diagram of the partially ordered set $(K, \subseteq)$. Here is a brief reminder of the basic facts and definitions of discrete Morse theory.

\medskip

Let $K$ be a simplicial complex. Its $p$-dimensional simplices
($p$-simplices for short) are denoted by $\alpha^p, \alpha^p_i, \beta^p,
\sigma^p$, etc. A \textit{discrete vector field} is a set of
pairs $D = \{\dots, (\alpha^p,\beta^{p+1}), \dots\}$ (called a matching) such
that:
\begin{enumerate}
 \setlength\itemsep{-1.5mm}
    \item[(a)]  each simplex of the complex participates in at most one
    pair;
    \item[(b)]  in each pair $(\alpha^p,\beta^{p+1})\in D$, the simplex $\alpha^p$ is a facet of $\beta^{p+1}$;
    \item[(c)]  the empty set $\emptyset\in K$ is not
matched, i.e.\ if $(\alpha^p,\beta^{p+1})\in D$ then $p\geq 0$.
\end{enumerate}
\noindent The pair $(\alpha^p, \beta^{p+1})$ can be informally
thought of as a vector in the vector field $D$. For this reason it
is occasionally denoted by $\alpha^p \rightarrow \beta^{p+1}$ or $\alpha^p \nearrow \beta^{p+1}$  (and
in this case $\alpha^p$ and $\beta^{p+1}$ are informally referred to as  the {\em beginning} and {\em the end}\/ of
the arrow $\alpha^p \rightarrow \beta^{p+1}$).

\smallskip
Given a discrete vector field $D$, a \textit{gradient path} in $D$
is a sequence of simplices (a zig-zag path)
$$\alpha_0^p \nearrow \beta_0^{p+1} \searrow \alpha_1^p \nearrow \beta_1^{p+1} \searrow \alpha_2^p \nearrow \beta_2^{p+1} \searrow \,  \cdots \,   \searrow\alpha_m^p \nearrow \beta_m^{p+1} \searrow \alpha_{m+1}^p$$
satisfying the following conditions:
\begin{enumerate}\setlength\itemsep{-1.5mm}
\item  $\big(\alpha_i^p,\ \beta_i^{p+1}\big)$ is a pair in $D$ for each $i$;
\item for each $i = 0,\dots, m$ the simplex $\alpha_{i+1}^p$
 is a facet of $\beta_i^{p+1}$;
 \item for each $i = 0,\dots, m-1,
 \, \alpha_i\neq \alpha_{i+1}$.
\end{enumerate}

A path is \textit{closed} if $\alpha_{m+1}^p=\alpha_{0}^p$. A
\textit{discrete Morse function } (DMF for short) is a discrete
vector field without closed paths.

\medskip
Assuming that a discrete Morse function is fixed, the {\em
critical simplices} are those simplices of the complex that are
not matched. The Morse inequality \cite{Forman2} implies that critical
simplices cannot be completely avoided.

A discrete Morse function $D$ is \textit{perfect}
if the number of critical $k$-simplices  equals the $k$-th
Betty number of the complex. It follows that $D$ is a perfect Morse function if and only if the number of all critical simplices equals the sum of all Betty numbers of $K$.

\medskip
A central idea of discrete Morse theory, as summarized in
the following theorem of R.~Forman, is to contract all matched
pairs of simplices and to reduce the simplicial complex $K$ to a
cell  complex (where critical simplices turn to the cells).

\begin{thm}\cite{Forman2} \label{ThmWedge} Assume that a
discrete Morse function on a simplicial complex $K$ has a single
zero-dimensional critical simplex $\sigma^0$ and that all other
critical simplices have the same dimension $N>1$. Then $K$ is
homotopy equivalent to a wedge of $N$-dimensional spheres.

More generally, if all critical simplices, aside from $\sigma^0$, have
dimension $\geq N$, then the complex $K$ is $(N-1)$-connected.  \qed
\end{thm}

\section*{Appendix 2. Comparison principle for equivariant maps }
\label{sec:Degrees}

The following theorem is proved in \cite{KB} (Theorem 2.1 in Section 2). Note that the condition that the $H_i$-fixed point sets $S^{H_i}$ are {\em locally $k$-connected} for $k\leq {\rm dim}(M^{H_i})-1$ is automatically satisfied if $S$ is a representation sphere. So in this case it is sufficient to show that the sphere $S^{H_i}$ is (globally) $({\rm dim}(M^{H_i})-1)$-connected which is equivalent to the condition
$${\rm dim}(M^{H_i}) \leq {\rm dim}(S^{H_i})\, (i=1,\dots, m) \, .$$

\begin{thm}\label{thm:K-B-2.1}
Let $G$ be a finite group acting on a compact topological manifold $M = M^n$ and on a sphere $S \cong S^n$ of the same dimension.
Let $N\subset M$ be a closed invariant subset and let $(H_1), (H_2), \dots, (H_k)$ be the orbit types in $M\setminus N$. Assume that the set $S^{H_i}$ is both globally and locally $k$-connected for all $k=0, 1,\dots, {\rm dim}(M^{H_i})-1$, where $i = 1,\dots, k$.
Then for every pair of $G$-equivariant maps $\Phi, \Psi : M\longrightarrow S$, which are equivariantly homotopic on $N$, there is the following relation
\begin{equation}\label{eqn:fundam-congruence}
 {\rm deg}(\Psi) \equiv {\rm deg}(\Phi)  \quad ({\rm mod}\, GCD\{\vert G/H_1\vert, \dots, \vert G/H_k\vert\}) \, .
\end{equation}
\end{thm}

\end{document}